\documentstyle[amscd,amssymb,verbatim,12pt,diagrams]{amsart}
\newcommand{\ad}{\operatorname{ad}}
\newcommand{\tr}{\operatorname{tr}}
\newcommand{\bfE}{{\bf E}}
\newcommand{\Div}{\operatorname{Div}}
\newcommand{\Aut}{\operatorname{Aut}}

\renewcommand{\mod}{\operatorname{mod}}

\newcommand{\coker}{\operatorname{coker}}

\newcommand{\SL}{\operatorname{SL}}

\newcommand{\lan}{\langle}
\newcommand{\ran}{\rangle}

\newcommand{\CC}{{\cal C}}

\newcommand{\supp}{\operatorname{supp}}

\newcommand{\Mat}{\operatorname{Mat}}

\newcommand{\Th}{\Theta}
\newcommand{\si}{\sigma}

\newcommand{\de}{\delta}
\newcommand{\eps}{\epsilon}

\renewcommand{\ker}{\operatorname{ker}}

\numberwithin{equation}{section}

\newtheorem{thm}{Theorem}[section]
\newtheorem{prop}[thm]{Proposition}
\newtheorem{lem}[thm]{Lemma}
\newtheorem{cor}[thm]{Corollary}
\newenvironment{rem}{\vspace{3mm}\noindent
{\bf Remark.}}{\vspace{3mm}}

\newcommand{\Pf}{\noindent {\it Proof}}

\newcommand{\ov}{\overline}

\renewcommand{\Im}{\operatorname{Im}}

\newcommand{\rk}{\operatorname{rk}}

\newcommand{\FF}{{\cal F}}
\newcommand{\HH}{{\cal H}}

\newcommand{\LL}{{\cal L}}
\newcommand{\ZZ}{{\cal Z}}
\newcommand{\UU}{{\cal U}}
\newcommand{\Om}{\Omega}

\newcommand{\Hom}{\operatorname{Hom}}
\newcommand{\Ext}{\operatorname{Ext}}
\newcommand{\End}{\operatorname{End}}

\renewcommand{\a}{\alpha}

\newcommand{\la}{\lambda}
\newcommand{\th}{\theta}
\newcommand{\C}{{\Bbb C}}
\newcommand{\R}{{\Bbb R}}
\newcommand{\Z}{{\Bbb Z}}
\newcommand{\Q}{{\Bbb Q}}

\newcommand{\wt}{\widetilde}

\newcommand{\sub}{\subset}
\newcommand{\ed}{\qed\vspace{3mm}}
\newcommand{\bfe}{{\bf e}}
\newcommand{\bfs}{{\bf s}}
\newcommand{\Exp}{\operatorname{Exp}}
\newcommand{\ovv}{\overline{v}}
\newcommand{\pa}{\partial}

\title{Analogues of the exponential map 
associated with complex structures on noncommutative two-tori}
\author{A. Polishchuk}
\thanks{Supported in part by NSF grant}

\begin{document}
\begin{abstract} We define and study analogues of exponentials
for functions on noncommutative two-tori that depend on a choice
of a complex structure. The major difference with the commutative
case is that our noncommutative exponentials can be defined only for
sufficiently small functions. We show that
this phenomenon is related to the existence
of certain discriminant hypersurfaces in an irrational rotation algebra.
As an application of our methods we give a very explicit characterization
of connected components in the group of invertible elements of this
algebra.    
\end{abstract}
\maketitle

\bigskip

\centerline{\sc Introduction}

\medskip

In this paper we study some natural constructions for
functions on noncommutative two-tori equipped with a complex structure.
Recall that for every number $\th\in\R\setminus\Q$
the algebra $A_{\th}$ of smooth functions on the noncommutative
torus $T_{\th}$ (also known as {\it irrational rotation algebra})
consists of expressions $\sum_{(m,n)\in\Z^2}a_{m,n}U_1^mU_2^n$
where coefficients $a_{m,n}\in\C$ decrease rapidly at infinity
and the multiplication is performed using the rule 
$$U_1U_2=\exp(2\pi i\th) U_2U_1.$$
Given an element $\tau\in\C\setminus\R$ (following \cite{PS} we will
always assume that $\Im(\tau)<0$) we define a derivation
$$\de_{\tau}:A_{\th}\to A_{\th}:\sum a_{m,n}U_1^m U_2^n\mapsto
2\pi i\cdot\sum_{m,n} (m\tau+n)a_{m,n} U_1^m U_2^n$$
We consider $\de_{\tau}$ as a complex structure on $T_{\th}$ and
denote the resulting complex noncommutative torus by $T_{\th,\tau}$.

The main object of our study is the equation
\begin{equation}\label{maineqintr}
\de_{\tau}(x)=xa
\end{equation}
for $x\in A_{\th}$, where $a\in A_{\th}$ is given.
In the commutative case this equation is clearly related with
the exponential map on smooth functions. It turns out that there is
a local analogue of this map for $A_{\th}$. However, in noncommutative
case there seem to be serious reasons why the exponential map does not 
extend to all functions. For example, we show that (\ref{maineqintr})
has a nonzero solution iff $\tr(a)\in 2\pi i(\Z+\Z\tau)$, where
$\tr(\sum a_{m,n}U_1^mU_2^n)=a_{0,0}$ 
(see Corollary \ref{maineqcor}), but
these solutions are not necessarily invertible in $A_{\th}$ (as in 
the commutative case). 

The study of equation (\ref{maineqintr}) turns out to
be closely related to the study of holomorphic structures on
the trivial holomorphic bundle over $T_{\th,\tau}$. By
a holomorphic bundle on $T_{\th,\tau}$ we mean a right projective
module $E$ over $A_{\th}$ equipped with a $\de_{\tau}$-connection, i.e.,
a linear map $\ov{\nabla}:E\to E$ satisfying the Leibnitz rule
$\ov{\nabla}(ea)=\ov{\nabla}(e)a+e\de_{\tau}(e)$.
The category $\CC$ of holomorphic bundles on $T_{\th,\tau}$ was studied in
\cite{PS} and \cite{P}. The main result of \cite{P} identifies $\CC$
with the heart of a certain $t$-structure on the derived category of
coherent sheaves on the elliptic curve $\C/(\Z+\Z\tau)$.
This leads to a classification of holomorphic bundles on
$T_{\th,\tau}$ up to isomorphism. 
We apply this classification to the study of solutions of (\ref{maineqintr}),
in particular, to the question of existence of a solution $x\in A_{\th}^*$,
where $A_{\th}^*\sub A_{\th}$ is the set of invertible elements.
The main result here is Theorem \ref{logthm} stating that
the map 
$$A_{\th}^*/\C^*\to A_{\th}:x\mapsto x^{-1}\de_{\tau}(x)$$
identifies $A_{\th}^*/\C^*$ with $\Om_{\tau}+2\pi i(\Z+\Z\tau)$
where $\Om_{\tau}$ is a dense open subset in the hyperplane
$H=\{a\in A_{\th}:\ \tr(a)=0\}$. 
More precisely, we prove that $\Om_{\tau}$ 
is the complement to the {\it discriminant hypersurface}
$\Th_{\tau}\sub H$ consisting of $a\in H$ such that the equation
$$\de_{\tau}(x)+ax-xa=0$$
has a nontrivial solution $x\in H$. We call it a hypersurface since
it coincides with the zero locus of a global holomorphic section of
a holomorphic line bundle over $A_{\th}$ (induced by Quillen's determinant
line bundle on the space of Fredholm operators of index zero).
We study the structure of $\Th_{\tau}$ in more detail. In particular,
we present some results that make us believe
that $\Th_{\tau}$ should be irreducible in an
appropriate sense.

As a byproduct of our study we deduce the following statement
about the structure of the group $A_{\th}^*$. It is well known
that the group of connected components of $A_{\th}^*$ can be identified
with $\Z^2$. We prove that for every $x\in A_{\th}^*$ one has
$\tr(x^{-1}\de_{\tau}(x))\in 2\pi i(\Z+\Z\tau)\sub\C$ and that the
map $x\mapsto\tr(x^{-1}\de_{\tau}(x))$ induces an isomorphism
on the groups of connected components.

The paper is organized as follows. In section \ref{formalexp-sec}
we study formal exponential maps for $A_{\th}$.
First, we present the construction of a map
$$\bfE_l(\tau,\cdot):A_{\th}[[t]]/\C[[t]]\to 1+t A_{\th}[[t]]$$
that specializes to
$a\mapsto \exp(ta)/\tr(\exp(ta))$ in the commutative case.
Then in section \ref{norm-exp-sec} we construct and study the normalized
map
$$\Exp_l(\tau,\cdot):A_{\th}[[t]]\to 1+t A_{\th}[[t]]$$
that specializes to $\exp(ta)$ in the commutative case.
In section \ref{convexp-sec} we prove that our formal exponentials
converge in a neighborhood of zero and then study 
logarithmic derivatives of invertible elements and the discriminant
hypersurface. 

\noindent
{\it Acknowledgment.} I am grateful to Chris Phillips for many
useful conversations and to Hanfeng Li for correcting a mistake in the
first draft.

\section{Formal exponentials}\label{formalexp-sec}

\subsection{Notation}

For every $v=(m,n)\in\Z^2$ we set 
$$U_v=\exp(-\pi i\th mn)U_1^m U_2^n$$
We have the following product rule in $A_{\th}$:
$$U_v\cdot U_{v'}=\exp(2\pi i\lan v,v'\ran)U_{v+v'},$$
where for $v'=(m',n')$ we set
$$\lan v,v'\ran=\lan v,v'\ran_\th=\frac{1}{2}\th (mn'-m'n).$$
For an element $a=\sum_{v\in\Z^2}a_v U_v\in A_{\th}$ we define the
{\it support} of $a$ as the set of all $v\in\Z^2$ such that $a_v\neq 0$.
Note that elements supported on a fixed rank-$1$ subgroup
$\Z\sub\Z^2$ form a commutative subalgebra in $A_{\th}$.

Recall that $A_{\th}$ is equipped with the $\C$-antilinear
antiinvolution $*$ defined by $(U_1)^*=U_1^{-1}$, $(U_2)^*=U_2^{-1}$.
For $v\in\Z^2$ we have $(U_v)^*=U_{-v}$,
hence for $a=\sum_{v\in\Z^2} a_v U_v$ we get 
$a^*=\sum_v a_v^* U_v$, where $a_v^*=\ov{a_{-v}}$.

Let $\iota:\Z^2\to\C$ be the homomorphism sending $(m,n)$ to
$2\pi i(m\tau+n)$. Then for every $v\in\Z^2$ we have
$$\de_{\tau}(U_v)=\iota(v)\cdot U_v.$$

\subsection{First construction}

Our formal exponentials will live in the ring $A_{\th}[[t]]$,
where $t$ is a formal variable commuting with $A_{\th}$.
We extend $\de_{\tau}$ to a $\C[[t]]$-linear derivation of
$A_{\th}[[t]]$ and the trace to a $\C[[t]]$-linear functional
$\tr:A_{\th}[[t]]\to\C[[t]]$.

\begin{thm}\label{formalthm} 
For every $\tau\in\C\setminus\R$ there exists a unique map
$$\bfE_l(\tau,\cdot):A_{\th}[[t]]/\C[[t]]\to 1+tA_{\th}[[t]],$$
with the following properties:

\noindent
(i) $\de_{\tau}(\bfE_l(\tau,a))=t\bfE_l(\tau,a)\cdot \de_{\tau}(a)$;

\noindent
(ii) $\tr(\bfE_l(\tau,a))=1$.
\end{thm}
 
We need the following purely algebraic statement.

\begin{lem}\label{alglem} 
Let $A$ be an associative algebra over $\Q$,
$M$ an $A$-bimodule, and $d:A\to M$ a derivation.
Assume that $a_1,\ldots,a_n$ are elements of $A$ such that
$$d(a_i)=a_{i-1}d(a_1) \text{ for } i=2,\ldots,n.$$
Then 
$$a_nd(a_1)\in [A,M]+d(A),$$
where $[A,M]\sub M$ is the linear span of the elements of the form
$[a,m]=am-ma$ for $a\in A$, $m\in M$.
\end{lem}

\Pf . For every $n_1,\ldots,n_k\in [1,n]$ let us denote
$$[n_1,\ldots,n_k]=a_{n_1}\ldots a_{n_{k-1}}d(a_{n_k})\mod [A,M]+d(A)$$
From the Leibnitz identity we get
$$
[n_1,\ldots,n_k]+c.p.(1,\ldots,k)=0
$$
in $A/([A,M]+d(A))$, where $c.p.(1,\ldots,k)$ denotes the terms
obtained by cyclic permutation from $(1,\ldots,k)$.
Since $k$ is invertible in $A$, this implies that
\begin{equation}\label{totsumeq}
\sum_{n_1+\ldots+n_k=n}[n_1,\ldots,n_k]=0.
\end{equation}
On the other hand, the identity $d(a_i)=a_{i-1}d(a_1)$ implies that
\begin{equation}\label{indeq}
[n_1,\ldots,n_k]=[n_1,\ldots,n_{k-1},n_k-1,1]
\end{equation}
for $n_k>1$. Let us denote
$$b_k=\sum_{n_1+\ldots+n_k=n}[n_1,\ldots,n_k,1].$$
For example, $b_1=[n,1]=a_n d(a_1)$, $b_{n+1}=0$.
We claim that for every $k$
$$b_k=-b_{k-1}.$$
Indeed, using (\ref{indeq}) we get
$$b_k=\sum_{n_1+\ldots+n_k=n}[n_1,\ldots,n_{k-1},n_k+1]=
\sum_{n_1+\ldots+n_k=n+1,n_k>1}[n_1,\ldots,n_k].$$
Now applying (\ref{totsumeq}) we can replace this with
$$b_k=-\sum_{n_1+\ldots+n_k=n+1,n_k=1}[n_1,\ldots,n_k]=-b_{k-1}$$
which proves out claim.
Therefore, 
$$a_n d(a_1)=b_1=\pm b_{n+1}=0$$ 
modulo $[A,M]+d(M)$.
\ed

\noindent
{\it Proof of Theorem \ref{formalthm}.}
Let us first show the uniqueness of $\bfE_l(\tau,a)$.
If $x,y\in 1+tA_{\th}[[t]]$ satisfy $x^{-1}\de(x)=y^{-1}\de(y)=t\de(a)$
then
$$\de(xy^{-1})=\de(x)y^{-1}-xy^{-1}\de(y)y^{-1}=0$$ 
hence $x=cy$ with $c\in\C[[t]]$. Furthermore,
$$1=\tr(x)=\tr(cy)=c\tr(y)=c,$$
so $x=y$.
To show the existence let us write
$\bfE_l(\tau,a)=1+a_1t+a_2t^2+\ldots$ with some $a_i\in A_{\th}[[t]]$, such that
$\tr(a_i)=0$ (so that condition (ii) is satisfied) and such that
$$\de_{\tau}(a_i)=a_{i-1}\de_{\tau}(a) \text{ for } i\ge 1,$$
where $a_0=1$. It is clear that condition (i) would follow, so it is enough to
show the existence of such $a_i$'s. 
To this end we set $a_1=a-\tr(a)$ and apply Lemma \ref{alglem}
to construct $a_i$ for $i\ge 2$ inductively. 
Note that we have  
$$[A_{\th},A_{\th}]\sub H=\de_{\tau}(A_{\th})$$
where $H\sub A_{\th}$ is the space of elements $x$ with $\tr(x)=0$.
A similar inclusion holds for the ring of formal power series.
Thus, when $a_1,\ldots,a_n$ are already constructed,
Lemma \ref{alglem} allows us to conclude that
$a_n\de_{\tau}(a)\in \de_{\tau}(A_{\th}[[t]])$, so we can define
$a_{n+1}$ as the unique element satisfying
$\de_{\tau}(a_{n+1})=a_n\de_{\tau}(a)$ and $\tr(a_{n+1})=0$.
\ed

\begin{rem} In the commutative case we have
$$\bfE_l(\tau,a)=\exp(ta)/\tr(\exp(ta))$$
due to property (ii). 
In section \ref{norm-exp-sec} we will present a
way of normalizing noncommutative exponentials which
reduces to the standard normalization in the commutative case.
\end{rem}

Similarly, one can define the function
$\bfE_r(\tau,\cdot)$ such that
$$\de_{\tau}(\bfE_r(\tau,a))=t\de_{\tau}(a)\cdot\bfE_r(\tau,a)$$
and $\tr(\bfE_r(\tau,a))=1$.
Let us extend the standard antiinvolution $*:A_{\th}\to A_{\th}$ to
an antiinvolution of $A_{\th}[[t]]$ by $t$-linearity.
From the equation $\de_{\tau}(a)^*=\de_{\ov{\tau}}(a^*)$ one can easily
derive that
\begin{equation}\label{bfE*-eq}
\bfE_r(\tau,a)^*=\bfE_l(\ov{\tau},a^*).
\end{equation}

\begin{prop}\label{prodprop} 
For every $a\in A_{\th}[[t]]$ one has
$$\bfE_l(\tau,a)\cdot\bfE_r(\tau,-a)=
\bfE_r(\tau,-a)\cdot\bfE_l(\tau,a)=:\bfs(\tau,a)\in 1+t\C[[t]].$$
\end{prop}

\Pf . Indeed, let $x=\bfE_l(\tau,a)$, $y=\bfE_r(\tau,-a)$, $\de=\de_{\tau}$. 
Then
$$\de(xy)=\de(x)y+x\de(y)=tx\de(a)y-tx\de(a)y=0.$$
This implies that $z=xy\in 1+t\C[[t]]$. 
Since $x$ is invertible in $A_{\th}[[t]]$ and
$z$ commutes with $x^{-1}$, we derive that $y=zx^{-1}$ and so
$yx=z=xy$. 
\ed

\begin{cor} One has $\ov{\bfs(\tau,a)}=\bfs(\ov{\tau},-a^*)$.
\end{cor}

\Pf . We have 
$$\ov{\bfs(\tau,a)}=\bfE_r(\tau,-a)^*\bfE_l(\tau,a)^*=\bfE_l(\ov{\tau},-a^*)
\bfE_r(\ov{\tau},a^*)=\bfs(\ov{\tau},-a^*).$$
\ed

\subsection{Formula in terms of coefficients}

Below we are going to write an explicit formula for $\bfE_l(\tau,a)$ in 
terms of the coefficients of $a$.

Recall that $\iota:\Z^2\to\C$ is the homomorphism sending $(m,n)$ to
$2\pi i(m\tau+n)$.
Let us define a collection of symmetric functions $f_n(v_1,\ldots,v_n)$
of $n$ lattice vectors recursively by setting $f_0=f_1(v)\equiv 1$,
$$f_n(v_1,\ldots,v_n)=
\iota(v_1+\ldots+v_n)^{-1}\sum_{i=1}^n\iota(v_i)f_{n-
1}(v_1,\ldots,\widehat{v_i},\ldots,v_n)
\exp(2\pi i\lan \sum_{j\neq i} v_j ,v_i\ran)
$$
if $v_1+\ldots+v_n\neq 0$, and $f_n(v_1,\ldots,v_n)=0$ otherwise. 
It is easy to see that $f_n$ can also be given by the following formula:
$$f_n(v_1,\ldots,v_n)=\sum_{\si\in S(v_1,\ldots,v_n)}
\frac{\iota(v_{\si(1)})\ldots \iota(v_{\si(n)})
\exp(2\pi i\sum_{i<j}\lan v_{\si(i)},v_{\si(j)}\ran)}
{\iota(v_{\si(1)})\iota(v_{\si(1)}+v_{\si(2)})\ldots
\iota(v_{\si(1)}+\ldots+v_{\si(n)})},$$
where the summation is taken over the set $S(v_1,\ldots,v_n)\sub S_n$ of all
permutations $\si$ for which the denominator 
in the corresponding term does not vanish.
When we want to stress the dependence of $f_n$ on $\th$ we write
$f_n(v_1,\ldots,v_n;\th)$ instead of $f_n(v_1,\ldots,v_n)$.

Let $\Div(\Z^2)_{\ge 0}$ denote the 
semigroup of effective divisors on $\Z^2$, i.e.,
of formal linear combinations $D=n_1(v_1)+\ldots+n_k(v_k)$ with $n_i\ge 0$ and
$v_i\in\Z^2$. We set 
$$\supp(D)=\{v_i: n_i>0\}$$
$$\deg(D)=n_1+\ldots+n_k,$$ 
$$s(D)=n_1v_1+\ldots+n_kv_k\in\Z^2,$$
$$D!=n_1!\ldots n_k!$$

\begin{thm}\label{formal2thm} 
For an element $a=\sum_{v\neq 0} a_v U_v\in A_{\th}[[t]]$, such that 
$a_0=0$, one 
has
$$\bfE_l(\tau,a)=
\sum_{D\in\Div(\Z^2)_{\ge 0}}\frac{t^{\deg(D)}}{D!}c(D)a_D U_{s(D)},$$
where for $D=(v_1)+\ldots+(v_d)$ (with not necessarily distinct $v_i$) we set
$$c(D)=c_{\th}(D):=f_d(v_1,\ldots,v_d;\th),$$
$$a_D=a_{v_1}\ldots a_{v_d}.$$
\end{thm}

\Pf . Recall that $\bfE_l(\tau,a)=\sum a_i t^i$, where
$\de(a_i)=a_{i-1}\de(a)$ and $\tr(a_i)=0$ for $i\ge 1$. 
This easily implies that
$$a_n=\sum_{D\in\Div(\Z^2)_{\ge 0}:\deg D=n}\frac{c(D)}{D!} a_D U_{s(D)}$$
for some constants $c(D)$ that do not depend on $a$. Note that the $c(D)$
are uniquely determined from this equation since the functions
$a\mapsto a_D$ for $D\in\Div(\Z^2)_{\ge 0}$ are linearly independent.
Note also that we have $c(D)=0$ if $s(D)=0$.
Comparing the coefficients with $a_D U_{s(D)}$ in the equation $\de(a_d)=a_{d-
1}\de(a)$,
where $d=\deg D$, we obtain the equation
$$\iota(s(D))\frac{c(D)}{D!}=\sum_{v\in\supp(D)}\iota(v)\frac{c(D-v)}{(D-v)!}
\exp(2\pi i\lan s(D-v),v\ran).$$
Equivalently,
\begin{equation}\label{cDeq}
\iota(s(D))c(D)=\sum_{v\in\supp(D)}n_v\iota(v)c(D-v)\exp(2\pi i\lan 
s(D-v),v\ran),
\end{equation}
where $n_v$ is the multiplicity of $v$ in $D$.
This implies that
$c(D)=f_d(v_1,\ldots,v_d)$ for $D=(v_1)+\ldots+(v_d)$.
\ed

\begin{cor}\label{fncor} For every $v_1,\ldots, v_n\in\Z^2$ such that
$v_1+\ldots+v_n=0$ one has
$$\sum_{i=1}^n\iota(v_i)f_{n-1}(v_1,\ldots,\widehat{v_i},\ldots,v_n)=0.$$
\end{cor}

\Pf . Apply equation (\ref{cDeq}) to $D=(v_1)+\ldots+(v_n)$.
\ed

\begin{rem} 
We do not know a direct combinatorial proof of 
the identity of Corollary \ref{fncor} except in the case when
$v_1,\ldots,v_n$ are sufficiently generic. See the remark after
Theorem \ref{formal-norm-thm}.
\end{rem}

Recall that $(\sum_v a_v U_v)^*=\sum_v \ov{a_{-v}} U_v$.
Therefore, using \eqref{bfE*-eq}
and Theorem \ref{formal2thm} we get
$$\bfE_r(\tau,a)=
\sum_{D\in\Div(\Z^2)_{\ge 0}}\frac{t^{\deg(D)}}{D!}c_{-\th}(D)a_D U_{s(D)}.$$

\subsection{Normalized exponential}\label{norm-exp-sec}

In this section we define a natural modification of $\bfE_l(\tau,\cdot)$
that reduces to the usual exponential in the commutative case.

The idea is to modify the symmetric functions $f_n(v_1,\ldots,v_n)$ that appear
in the explicit formula for $\bfE_l(\tau,\cdot)$.
Let us consider the rational functions
$$R_n(x_1,\ldots,x_n)=\frac{x_1x_2\ldots x_n}
{x_1(x_1+x_2)\ldots(x_1+\ldots+x_n)}.$$
For a permutation $\si\in S_n$ and an $n$-tuple $x=(x_1,\ldots,x_n)$
let us denote $x^{\si}=(x_{\si(1)},\ldots,x_{\si(n)})$.
Recall that for an $n$-tuple of lattice vectors $\ovv=(v_1,\ldots,v_n)$
we have
$$f_n(\ovv)=\sum_{\si\in S(\ovv)}
R_n(\iota(\ovv)^{\si})\exp(2\pi i\sum_{i<j}
\lan v_{\si(i)},v_{\si(j)}\ran)$$
where $\iota(\ovv)=(\iota(v_1),\ldots,\iota(v_n))$, and
$S(\ovv)\sub S_n$ is the set of all permutations $\si$
such that $R_n$ is regular at $\ovv^{\si}$.
We want to define a modified expression $f^*_n(\ovv)$ 
that takes into account contributions from $\si$ such
that $R_n$ is not defined at $\ovv^{\si}$. For this we will combine
several terms corresponding to such permutations, so that the poles will
cancel out.

An $n$-tuple $\ovv=(v_1,\ldots,v_n)$ determines a map
$$t_{\ovv}:S_n\to\Z^2\otimes_{\Z}\Z^2:
\si\mapsto\sum_{i<j}v_{\si(i)}\otimes v_{\si(j)}.$$
The fibers of $t_{\ovv}$ define a partition
$$S_n=S(\ovv,1)\sqcup \ldots\sqcup S(\ovv,r)$$
of $S_n$. 
The crucial observation is
the following remarkable cancellation of poles.

\begin{prop}\label{nopoleprop} The rational functions 
$$R_{S(\ovv,j)}(x):=\sum_{\si\in S(\ovv,j)}R_n(x^{\si})$$
for $j=1,\ldots,r$ are regular at $\iota(\ovv)$. More precisely,
for every collection $1\le i_1<\ldots<i_r\le n$ such that
$v_{i_1}+\ldots+v_{i_r}=0$ the functions $R_{S(\ovv,j)}$
do not have poles along the hyperplane $x_{i_1}+\ldots+x_{i_r}=0$.
\end{prop}

The proof is based on two lemmata.

\begin{lem}\label{cyclic-inv-lem} 
If $v_1+\ldots+v_n=0$ then
$$\sum_{i<j} v_i\otimes v_j=\sum_{i<j} v_{i+1}\otimes v_{j+1},$$
where we view indices as elements of $\Z/n\Z$, so that
$v_{n+1}=v_1$.
\end{lem}

The proof is straightforward.

\begin{lem}\label{cyclic-reg-lem} 
Let $\zeta\in S_n$ be the cyclic permutation of order $n$.
Then the function
$$\sum_{i=0}^{n-1} R_n(x^{\zeta^i})$$
has no pole along $x_1+\ldots+x_n=0$. 
\end{lem}

\Pf . As before we view indices as elements of $\Z/n\Z$.
It suffices to prove that the function
$$\sum_{i=1}^n\frac{1}{x_i(x_i+x_{i+1})\ldots (x_i+x_{i+1}+\ldots+x_{i+n-1})}$$
has no pole along $x_1+\ldots+x_n=0$. Since
$$x_i+x_{i+1}+\ldots+x_j\equiv -(x_{j+1}+x_{j+2}+\ldots+x_{i-1})\mod
(x_1+\ldots+x_n)$$
we can replace the above expression with
$$\frac{1}{x_1+\ldots+x_n}\cdot\sum_{i=1}^n\frac{(-1)^{i-1}}
{(x_1+\ldots+x_{i-1})(x_2+\ldots+x_{i-1})x_{i-1}
x_i(x_i+x_{i+1})\ldots(x_i+\ldots+x_{n-1})}.$$
Note that the new sum depends only on $(x_1,\ldots,x_{n-1})$.
Hence, our assertion is equivalent to the identity
$$\sum_{i=1}^n\frac{(-1)^{i-1}}
{(x_1+\ldots+x_{i-1})(x_2+\ldots+x_{i-1})x_{i-1}\cdot
x_i(x_i+x_{i+1})\ldots(x_i+\ldots+x_{n-1})}=0.$$
Since the left-hand side is homogeneous of degree $-(n-1)$, it
suffices to prove that it has no poles. 
The only possible poles are along the hyperplanes 
$h_{ij}=x_i+x_{i+1}+\ldots+x_j=0$.
Now one can check easily that for every $i\le j$ there
are exactly two terms in the above sum that have poles along $h_{ij}=0$
and that their polar parts cancel out.
\ed

\noindent
{\it Proof of Proposition \ref{nopoleprop}.}
We use induction in $n$ and the following recursive formula for $R_n$:
$$R_n(x_1,\ldots,x_n)=R_{n-1}(x_1,\ldots,x_{n-1})\frac{x_n}{x_1+\ldots+x_n}.$$
Assume that the assertion holds for $n-1$. We can write
$R_{S(\ovv,j)}$ in the form
$$R_{S(\ovv,j)}=\sum_{m=1}^n \sum_{\si\in S(\ovv,j):\si(n)=m}
R_{n-1}(x_{\si(1)},\ldots,x_{\si(n-1)})\frac{x_m}{x_1+\ldots+x_n}$$
For a given $m$ let us denote $S(\ovv,j;m)$ 
the set of $\si\in S(\ovv,j)$ such
that $\si(n)=m$. Fix $\si_0\in S(\ovv,j;m)$. Then
$S(\ovv,j;m)\si_0^{-1}$ leaves $m$ stable. Furthermore, it is clear
that $S(\ovv,j;m)\si_0^{-1}$ coincides with one piece 
in the partition of the set of
permutations of $(1,\ldots,\widehat{m},\ldots,n)$ associated
with the $(n-1)$-tuple $(v_1,\ldots,\widehat{v_m},\ldots,v_n)$.
Thus, by the induction assumption all the sums
$$\sum_{\si\in S(\ovv,j;m)}R_{n-1}(x_{\si(1)},\ldots,x_{\si(n-1)})$$
have no poles along the required hyperplanes. 
If $v_1+\ldots+v_n\neq 0$ this finishes the proof.
It remains to consider the case when $v_1+\ldots+v_n=0$.
By Lemma \ref{cyclic-inv-lem} in this case all the
sets $S(\ovv,i)$ are right cosets for the subgroup $\Z/n\Z\sub S_n$
generated by the cyclic permutation $\zeta\in S_n$, where 
$\zeta(i)=i+1$. Therefore, by Lemma \ref{cyclic-reg-lem} the functions
$R_{S(\ovv,j)}$ have no poles along the hyperplane $x_1+\ldots+x_n=0$.
\ed

Let us define a rational function of $x=(x_1,\ldots,x_n)$ by
$$R_{\ovv}(x)=\sum_{\si\in S_n}R(x^{\si})
\exp(2\pi i \lan t_{\ovv}(\si)\ran)=
\sum_{j=1}^r R_{S(\ovv,j)}(x)
\exp(2\pi i\lan t_{\ovv}(S(\ovv,j))\ran),$$
where $\lan t_{\ovv}(S(\ovv,j))\ran$ is the common value of
$$\lan t_{\ovv}(\cdot)\ran:=\lan \cdot,\cdot\ran\circ t_{\ovv}:
S_n\to\R$$ on $S(\ovv,j)$.
By Proposition \ref{nopoleprop} this function is regular at $\iota(\ovv)$,
so we can set
$$f^*_n(\ovv)=f^*_n(\ovv;\th)=R_{\ovv}(\iota(\ovv)).$$
It is easy to check that $f^*_n$ is symmetric in $v_1,\ldots,v_n$.
As before we have $f^*_0=f^*_1(v)\equiv 1$.

\begin{thm}\label{formal-norm-thm} 
For an element $a=\sum_v a_v U_v\in A_{\th}[[t]]$
consider the series
$$\Exp_l(\tau,a):=\sum_{D\in\Div(\Z^2)_{\ge 0}}
\frac{t^{\deg(D)}}{D!} c^*(D)a_D U_{s(D)}\in 1+t A_{\th}[[t]],$$
where for $D=(v_1)+\ldots+(v_d)$ we set
$$c^*(D)=c^*_{\th}(D)=f^*_n(v_1,\ldots,v_d;\th).$$
Then 

\noindent
(i) $\de_{\tau}(\Exp_l(\tau,a))=t\Exp_l(\tau,a)\cdot\de_{\tau}(a)$;

\noindent
(ii) if $\th=0$ or if the support of $a$ is contained in
a rank-$1$ subgroup $\Z\sub\Z^2$ then $\Exp_l(\tau,a)=\exp(ta)$;

\noindent
(iii) for $a\in A_{\th}[[t]]$ and $z\in\C[[t]]$ we have
\begin{equation}\label{exp-scalar-eq}
\Exp_l(z+a)=\exp(tz)\cdot\Exp_l(a).
\end{equation}
\end{thm}

\Pf . As in Theorem \ref{formal2thm} the proof of (i)
reduces to the following identity:
$$\iota(v_1+\ldots+v_n)f^*_n(v_1,\ldots,v_n)=
\sum_{m=1}^n\iota(v_m)f^*_{n-1}(v_1,\ldots,\widehat{v_m},\ldots,v_n)
\exp(2\pi i\lan\sum_{j\neq m} v_j,v_m\ran).$$
We claim that in fact, there is an identity between the corresponding
rational functions:
$$(x_1+\ldots+x_n)R_{\ovv}(x_1,\ldots,x_n)=
\sum_{m=1}^n x_m R_{\ovv^{(m)}}(x_1,\ldots,\widehat{x_m},\ldots,x_n)
\exp(2\pi i\lan\sum_{j\neq m} v_j,v_m\ran),$$
where $\ovv^{(m)}=(v_1,\ldots,\widehat{v_m},\ldots,v_n)$.
Indeed, we can rewrite the left-hand side as
\begin{align*}
&(x_1+\ldots+x_n)
\sum_{\si\in S_n}R_n(x^{\si})\exp(2\pi i\lan t_{\ovv}(\si)\ran)=\\
&\sum_{\si\in S_n}x_{\si(n)}
R_{n-1}(x_{\si(1)},\ldots,x_{\si(n-1)})
\exp(2\pi i[\sum_{i<j<n}\lan v_{\si(i)}, v_{\si(j)}\ran
+\lan\sum_{j\neq m} v_j,v_m\ran])=\\
&\sum_{m=1}^n x_m R_{\ovv^{(m)}}(x_1,\ldots,\widehat{x_m},\ldots,x_n)
\exp(2\pi i\lan\sum_{j\neq m} v_j,v_m\ran).
\end{align*}

If $\th=0$ or if vectors
$v_1,\ldots,v_n$ belong to a rank-$1$ subgroup $\Z\sub\Z^2$ then we have
$$R_{\ovv}(x)=\sum_{\si\in S_n}R_n(x^{\si}).$$
Applying Proposition \ref{nopoleprop} we deduce that this
function has no poles. Since it is homogeneous of degree $0$
it must be a constant. To compute this constant we substitute
$x_1=\ldots=x_n=1$ and obtain that in this case
$R_{\ovv}(x)\equiv 1$. This immediately implies (ii).

To prove part (iii) it suffices to check \eqref{exp-scalar-eq}
as a formal identity in $A_{\th}[[z,t]]$. We claim
that both sides solve the same differential
equation $\frac{\pa f}{\pa z}(z,t)=tf(z,t)$ with the same initial
condition $f(0,t)=\Exp_l(\tau,a)$. Indeed, this is clear for the
right-hand side of \eqref{exp-scalar-eq}, 
so we just have to check this for the left-hand side. 
It is easy to check that
$$\frac{1}{D!}\cdot\frac{\pa [(a+z)_D]}{\pa z}
=\frac{(a+z)_{D-(0)}}{(D-(0))!}$$
if $D-(0)\ge 0$, otherwise the derivative is zero,
This implies that the derivative of $\Exp_l(z+a)$ with respect to $z$
is equal to $t\Exp_l(z+a)$ as we claimed.
\ed

\begin{rem} For $n$ lattice vectors $(v_1,\ldots,v_n)$ such that
$v_1+\ldots+v_n=0$ we obtain from the above proof that
$$\sum_{m=1}^n\iota(v_m)f^*_{n-1}(v_1,\ldots,\widehat{v_m},\ldots,v_n)=0.$$
If $v_1,\ldots,v_n$ are sufficiently generic this reduces to the
identity of Corollary \ref{fncor}.
\end{rem}

\begin{cor}\label{exp-uniq-cor} 
If $\Exp_l(\tau,a)=\Exp_l(\tau,b)$ for some $a,b\in A_{\th}[[t]]$
then $a=b$.
\end{cor}

\Pf . Using part (i) of the above theorem we immediately deduce that
$b=a+z$ for some $z\in\C[[t]]$. Now part (iii) implies that
$\exp(tz)=1$ and hence $z=0$.
\ed

We can also define the right exponential map
$\Exp_r(\tau,\cdot)$ by setting
$$\Exp_r(\tau,a)=\Exp_l(\ov{\tau},a^*)^*.$$
Then 
$$\de_{\tau}(\Exp_r(\tau,a))=t\de_{\tau}(a)\cdot\Exp_r(\tau,a).$$
Using the definition of $\Exp_l(\tau,\cdot)$ we get
$$\Exp_r(\tau,a)=\sum_{D\in\Div(\Z^2)_{\ge 0}}
\frac{t^{\deg(D)}}{D!} c^*_{-\th}(D)a_D U_{s(D)}.$$

\begin{prop}\label{inv-prop} For every $a\in A_{\th}[[t]]$ one has
$$\Exp_l(\tau,a)^{-1}=\Exp_r(\tau,-a).$$
\end{prop}

\Pf . It suffices to prove that for every divisor
$D>0$ on $\Z^2$ one has
$$\sum_{D_1,D_2\in\Div(\Z^2)_{\ge 0}: D_1+D_2=D}
\frac{(-1)^{\deg D_2}}{D_1!D_2!}c^*_{\th}(D_1)c^*_{-\th}(D_2)
\exp(2\pi i\lan s(D_1),s(D_2)\ran_{\th})=0.$$
This is equivalent to the following identity that should
hold for every $n$-tuple $(v_1,\ldots,v_n)$ of lattice vectors
(where $n>0$):
\begin{align*}
&\sum_{I=\{i_1,\ldots,i_r\}\sub [1,n]}
(-1)^{n-r}f_r^*(v_{i_1},\ldots,v_{i_r};\th)
f_{n-r}^*(v_{j_1},\ldots,v_{j_{n-r}};-\th)\\
&\times
\exp(2\pi i\lan v_{i_1}+\ldots+v_{i_r},v_{j_1}+\ldots+v_{j_{n-r}}\ran_{\th})
=0,
\end{align*}
where the sum is taken over all subsets $I=\{i_1,\ldots,i_r\}$ 
of $[1,n]=\{1,2,\ldots,n\}$ and
we denote by $\{j_1,\ldots,j_{n-r}\}$ the complement to $I$.
Recalling the definition of 
the functions $f_n^*$ we see that it suffices to prove the following
identity between rational functions of $x_1,\ldots,x_n$:
\begin{align*}
&\sum_{r=0}^n\sum_{\si\in S_n}
(-1)^{n-r}R_r(x_{\si(1)},\ldots,x_{\si(r)})
R_{n-r}(x_{\si(n)},\ldots,x_{\si(r+1)})\\
&\exp(2\pi i
[\sum_{i<j\le r}\lan v_{\si(i)},v_{\si(j)}\ran_{\th}+
\sum_{i>j>r}\lan v_{\si(i)},v_{\si(j)}\ran_{-\th}+
\lan\sum_{i\le r} v_{\si(i)},\sum_{j>r} v_{\si(j)}\ran_{\th}])=0.
\end{align*}
Due to the skew-symmetric of $\lan\cdot,\cdot\ran_{\th}$ we have
$$\lan v,v'\ran_{-\th}=\lan v',v\ran_{\th}.$$
Using this observation we can rewrite our identity as
$$\sum_{r=0}^n\sum_{\si\in S_n}
(-1)^{n-r}R_r(x_{\si(1)},\ldots,x_{\si(r)})
R_{n-r}(x_{\si(n)},\ldots,x_{\si(r+1)})
\exp(2\pi i\sum_{i<j}\lan v_{\si(i)},v_{\si(j)}\ran_{\th})=0.$$
Thus, it is enough to prove that
$$\sum_{r=0}^n (-1)^{n-r} R_r(x_1,\ldots,x_r)R_{n-r}(x_n,\ldots,x_{r+1})=0.$$
This is equivalent to the following identity:
\begin{align*}
&\sum_{i=0}^{n-1}\frac{(-1)^i}
{x_1(x_1+x_2)\ldots(x_1+\ldots+x_{n-i})x_n(x_n+x_{n-1})\ldots
(x_n+\ldots+x_{n-i+1})}=\\
&\frac{(-1)^{n-1}}{x_n(x_n+x_{n-1})\ldots(x_n+\ldots+x_1)}.
\end{align*}
We claim that more generally for $0<j<n$ one has
\begin{align*}
&\sum_{i=0}^{j}\frac{(-1)^i}
{x_1(x_1+x_2)\ldots(x_1+\ldots+x_{n-i})x_n(x_n+x_{n-1})\ldots
(x_n+\ldots+x_{n-i+1})}=\\
&\frac{(-1)^j}
{x_1(x_1+x_2)\ldots(x_1+\ldots+x_{n-j-1})\cdot(x_1+\ldots+x_n)\cdot
x_n(x_n+x_{n-1})\ldots(x_n+\ldots+x_{n-j+1})}.
\end{align*}
This can be easily checked by induction in $j$.
\ed

\begin{cor} One has
$$\tr(\Exp_l(\tau,a))\cdot\tr(\Exp_r(\tau,-a))=\bfs(\tau,a)^{-1}.$$
\end{cor}

\Pf . This follows from the relation
$\Exp_l(\tau,a)=\tr(\Exp_l(\tau,a))\bfE_l(\tau,a)$
and the similar relation for right exponentials.
\ed

Finally, let us state the analogue of the equation 
$\exp(x)\exp(y)=\exp(x+y)$ for our exponentials.

\begin{prop} For every $a,b\in A_{\th}[[t]]$ there exists a unique
$\varphi(a,b)\in A_{\th}[[t]]$ such that
\begin{equation}\label{exp-prod-eq}
\Exp_l(a)\Exp_l(b)=\Exp_l(\varphi(a,b)+b)
\end{equation}
and
$$\de_{\tau}(\varphi(a,b))=\Exp_l(\tau,b)^{-1}\de_{\tau}(a)\Exp_l(\tau,b).$$
\end{prop}

\Pf . Since the element $\Exp_l(\tau,b)^{-1}\de_{\tau}(a)\Exp_l(\tau,b)\in
A_{\th}[[t]]$ has zero trace we can find some element
$\varphi'\in A_{\th}[[t]]$ such that
$$\de_{\tau}(\varphi')=\Exp_l(\tau,b)^{-1}\de_{\tau}(a)\Exp_l(\tau,b).$$
Now one can easily check that the left-hand side of \eqref{exp-prod-eq} is
a solution of the equation 
$$\de_{\tau}(x)=tx\de_{\tau}(\varphi'+b).$$
It follows that
$$\Exp_l(a)\Exp_l(b)=f(t)\Exp_l(\varphi'+b)$$
for some $f(t)\in 1+t\C[[t]]$.
Writing $f(t)$ in the form $\exp(tz)$ for $z\in\C[[t]]$ and
using part (iii) of Theorem \ref{formal-norm-thm} we conclude
that \eqref{exp-prod-eq} holds for $\varphi(a,b)=z+\varphi'$.
The uniqueness follows from Corollary \ref{exp-uniq-cor}.
\ed

\subsection{Compatibility with $\SL_2(\Z)$-action}

The group $\SL_2(\Z)$ acts on the algebra $A_{\th}$ by automorphisms.
Namely, to 
an element $g=\left(\matrix a & b \\ c & d\endmatrix\right)\in\SL_2(\Z)$
one associates an automorphism
$$\a_g:A_{\th}\to A_{\th}: U_v\mapsto U_{gv}.$$
These automorphisms are compatible with derivations $\de_{\tau}$ and
the modular action of $\SL_2(\Z)$ on $\tau$:
$$\a_g\de_{\tau}\a_g^{-1}=\de_{g\tau}$$
where $g\tau=(a\tau+b)/(c\tau+d)$.

The only way in which the standard basis of $\Z^2$ enters into
our formulae for $\bfE_l(\tau,a)$ and $\Exp_l(\tau,a)$ is
through the definition of the homomorphism $\iota=\iota_{\tau}:\Z^2\to\C$.
For an element $g\in\SL_2(\Z)$ we have
$$\iota_{g\tau}(v)=\iota_{\tau}(\sideset{^t}{}{g}v)\cdot j_{\tau}(g)^{-1},$$
where $\sideset{^t}{}{g}$ is the transpose of $g$, and
$j_{\tau}(g)=c\tau+d$. Since the functions $f^*_n(v_1,\ldots,v_n)$
(resp., $f_n(v_1,\ldots,v_n)$) are homogeneous of degree $0$ this implies
that
$$f^*_n(v_1,\ldots,v_n;g\tau)=f^*_n(\sideset{^t}{}{g}v_1,\ldots,
\sideset{^t}{}{g}v_n;\tau)$$
and similarly for the functions $f_n$. This easily implies the following
result.

\begin{prop}
For every $a\in A_{\th}[[t]]$ one has
$$\bfE_l(g'\tau,\a_g(a))=\a_g\bfE_l(\tau,a),$$
$$\bfE_r(g'\tau,\a_g(a))=\a_g\bfE_r(\tau,a),$$
$$\Exp_l(g'\tau,\a_g(a))=\a_g\Exp_l(\tau,a),$$
$$\Exp_r(g'\tau,\a_g(a))=\a_g\Exp_r(\tau,a),$$
$$\bfs(g'\tau,\a_g(a))=\bfs(\tau,a),$$
where $g'=\sideset{^t}{^{-1}}{g}$.
\end{prop}

\section{Convergent exponentials,
logarithmic derivatives and the discriminant hypersurface}
\label{convexp-sec}

\subsection{Norms} 

In this section we establish notation for
some of the norms we use on $A_{\th}$.

First of all we have the operator algebra norm $\|\cdot\|$ on $A_{\th}$:
$\|a\|$ is defined to be the norm of the operator of multiplication
by $a$ on the $L^2$-completion of $A_{\th}$ (it does not matter whether
one considers left or right multiplication). Recall that the $C^*$-algebra
$\ov{A}_{\th}$ is defined as the completion of $A_{\th}$ with respect
to this norm.

We will also consider the $L^2$-norm 
$$\|a\|_0=(\sum_{v\in\Z^2}|a_v|^2)^{1/2},$$
where $a=\sum_{v\in\Z^2} a_v U_v$. 
By the definition, for every $a,b\in A_{\th}$
we have
$$\|ab\|_0\le \|a\|\cdot \|b\|_0.$$

For every $s\ge 0$ we can also consider the norm
$$\|a\|_s^2=\sum_{i=0}^s \|\de_{\tau}^i a\|_0^2.$$
The completion of $A_{\th}$ with respect to this norm
is the {\it Sobolev space} $W_s$. We have a sequence of
embeddings $W_0\supset W_1\supset W_2\supset\ldots$
These spaces enjoy all
the usual properties of the Sobolev spaces (see \cite{P}).
In particular, it can be extended to a sequence of spaces
$(W_s)$, where $s\in\Z$, so that $W_s$ is dual to $W_{-s}$.
Also, the intersection of all $W_s$ is $A_{\th}$ and the topology
on $A_{\th}$ is the one determined by the collection of seminorms
coming from the $W_s$.

Finally, at one point we will use the norm
$$\|a\|_{l^1}:=\sum_{v\in\Z^2}|a_v|.$$
For $a,b\in A_{\th}$ we have
$$\|ab\|_{l^1}\le\sum_{v\in\Z^2}\|a_v U_v b\|_{l^1}=
\sum_{v\in\Z^2}|a_v|\cdot \|b\|_{l^1}=\|a\|_{l^1}\cdot \|b\|_{l^1}.$$
On the other hand,
$$\|a\|\le\sum_{v\in\Z^2}\|a_v U_v\|=\|a\|_{l^1}.$$

\subsection{Convergence criteria}

If one wants to solve the equation 
\begin{equation}\label{mainequation}
\de_{\tau}(x)=x\de_{\tau}(a)
\end{equation}
for given $a\in A_{\th}$ one can try to consider the series
\begin{equation}\label{mainseries}
x=\bfE_l(\tau,a)|_{t=1}=1+a_1+a_2+\ldots
\end{equation}
obtained by substituting $t=1$ into $\bfE_l(\tau,a)=1+ta_1+t^2a_2+\ldots$. 
This series does not always converge even for $\th=0$
(remember that in the commutative case it represents the function
$\exp(ta)/\tr(\exp(ta))$). One of our goals in this section is to show that
\eqref{mainseries} always converges for sufficiently small $a$.

\begin{rem}
In the commutative case the convergence improves when we switch from
$\bfE_l(\tau,a)$ to $\Exp_l(\tau,a)$ 
since the latter series becomes the usual exponential. However, this
does not seem to work in noncommutative case.
Indeed, as we will see in Proposition \ref{linesprop}(ii)
there exists $a\in A_{\th}$ for which 
\eqref{mainequation} has no invertible solutions.
By Proposition \ref{inv-prop} this means that either
$\Exp_l(\tau,a)|_{t=1}$ or $\Exp_r(\tau,-a)|_{t=2}$ diverges.
\end{rem}

The next lemma follows from Corollary 7.16 of \cite{Schwe}. Since
our particular case is much easier we present here a more direct proof 
due to Chris Phillips.

\begin{lem}\label{invertlem} 
If $x\in A_{\th}$ is invertible as an element of $\ov{A}_{\th}$
then $x\in A_{\th}^*$.
\end{lem}

\Pf . First, we claim that for every $a\in A_{\th}$ such that
$\|a\|<1$, the series
$1+a+a^2+\ldots$ converges to an element of $A_{\th}$. Indeed,
it suffices to show that it converges absolutely in every Sobolev
space $W_s$. But this follows from the fact that for every $k$ there exists
a constant $C$ (depending on $a$) such that 
$\|\de^k(a^n)\|\le C\| a\|^{n-k}$ for $n>k$.
Thus, the assertion is true when $\|x-1\|<1$. The general case can
be reduced to this by the following trick.
Let $x$ be an element of $A_{\th}$ such that there exists an inverse
$x^{-1}\in\ov{A}_{\th}$. Choose $y\in A_{\th}$ such 
$\|y-x^{-1}\|<\|x\|^{-1}$. Then the element $xy\in A_{\th}$ satisfies
$\|xy-1\|<1$. Hence, $xy\in A_{\th}^*$. It follows
that $x$ has a right inverse in $A_{\th}$. Similar argument shows
that $yx\in A_{\th}^*$ and hence $x$ has a left inverse. This implies
that $x\in A_{\th}^*$.
\ed

\begin{rem} In fact, it is also true that if $xy=1$ in $A_{\th}$
then $yx=1$. This follows from the existence of
a finite positive faithful trace on $A_{\th}$ (cf. \cite{Dav}, second
half of p.~101 and Exer. 6 and 9 of ch.~4).
Here is another argument using Rieffel's classification of finitely
generated projective $A_{\th}$-modules in \cite{R-canc}
(we assume that $\th$ is irrational).
Let $I\sub A_{\th}$ the set of all $a$ such that $xa=0$. Then
$I$ is a right ideal and we have $A_{\th}=I\oplus yA_{\th}$.
It follows that $yA_{\th}$ and $I$ are projective modules.
Furthermore, $yA_{\th}$ is isomorphic to
$A_{\th}$ hence $\rk yA_{\th}=1$. It follows that $\rk I=0$, so $I=0$. 
\end{rem}

\begin{thm}\label{convthm} Let 
$$d=\min_{(m,n)\in\Z^2\setminus\{(0,0)\}}|m\tau+n|.$$
Then for every $a\in A_{\th}$ with $\tr(a)=0$ such that
$\|\de_{\tau}(a)\|<2\pi d$, the series \eqref{mainseries}
converges to an element $\bfe_l(\tau,a)\in A_{\th}$ satisfying the equation
$$\de_{\tau}(\bfe_l(\tau,a))=\bfe_l(\tau,a)\de_{\tau}(a).$$
Furthermore, if 
$\|\de_{\tau}(a)\|_{l^1}<\pi d$ then $\bfe_l(\tau,a)\in A_{\th}^*$.
The map $a\mapsto\bfe_l(\tau,a)$ from
$U_{\tau}=\{a\in H:\ \|\de_{\tau}(a)\|_{l^1}<\pi d\}$
to $A_{\th}^*$ is continuous.
\end{thm}

\Pf . Recall that $a_i$'s are defined inductively by the conditions
$\de(a_i)=a_{i-1}\de(a)$, $\tr(a_i)=0$, where $\de=\de_{\tau}$.
It follows that
$$\|a_i\|_0\le (2\pi d)^{-1} \|\de(a_i)\|_0\le 
(2\pi d)^{-1}\|\de(a)\|\cdot\|a_{i-1}\|_0.$$
Hence, $\|a_n\|_0\le (\frac{\|\de(a)\|}{2\pi d})^n$ for every $n\ge 0$.
It follows that the series $1+a_1+a_2+\ldots$ converges to an
element $x$ in the $L^2$-completion of $A_{\th}$.
Furthermore, since
$$\|\de(a_n)\|_0\le \|\de(a)\|\cdot\|a_{n-1}\|_0\le 
\|\de(a)\|\cdot(\frac{\|\de(a)\|}{2\pi d})^{n-1},$$
we see that in fact $x$ belongs to the Sobolev space $W_1$ 
and satisfies the equation
$$\de(x)=x\de(a).$$
This implies that $x\in A_{\th}$. 

To prove the last assertion we first observe that
$\|a_n\|_{l^1}\le r^n$ for all $n$, where
$r=\|\de(a)\|_{l^1}/(2\pi d)$ (this is deduced in the same way as above).
Hence, 
$$\|a_1+a_2+\ldots\|\le \|a_1+a_2+\ldots\|_{l^1}\le \frac{r}{1-r}<1$$
if $r<1/2$. Therefore, $x$ is invertible as an element of $\ov{A}_{\th}$.
It remains to apply Lemma \ref{invertlem}.
\ed

The following result shows that our exponentials converge if
all the coefficients of $a$ belong to closed
halfplanes not containing zero.

\begin{prop}\label{suppconvprop}
Let $a=\sum_v a_v U_v\in A_{\th}$ be an element with $\tr(a)=0$.
Assume that there exists a homomorphism $h:\Z^2\to\R$ and a positive
constant $\eps$ such that $h>\eps$ on
$\supp(a)=\{v\in\Z^2: a_v\neq 0\}$. Then the above series $1+a_1+a_2+\ldots$
converges to an element of $A_{\th}^*$. 
\end{prop}

\Pf . Since the elements $a_i\in H$
are defined inductively by $a_1=a$ and 
$\de(a_i)=a_{i-1}\de(a)$, we deduce that $h>n\eps$ on $\supp(a_n)$.
Let $d_n=\min_{(m,n): h(m,n)>n\eps}|m\tau+n|$. Then we have
$$\|a_n\|_0\le (2\pi d_n)^{-1}\|\de(a_n)\|_0\le (2\pi d_n)^{-1}\|\de(a)\|
\cdot \|a_{n-1}\|_0.$$
Since $d_n$ grows linearly with $n$ this immediately implies convergence
of the series $\sum a_n$ in the $L^2$-norm. Arguing as in the proof
of Theorem \ref{convthm} we derive that $x=1+a_1+a_2+\ldots$ belongs
to $A_{\th}^*$.
\ed

\begin{rem} In the situation of the above proposition one has
$\Exp_l(\tau,a)|_{t=1}=\bfE_l(\tau,a)|_{t=1}$.
However, at present we do not know any criteria for convergence
of $\Exp_l(\tau,a)|_{t=1}$ similar to Theorem \ref{convthm}.
\end{rem}

\subsection{Discriminant hypersurface}\label{hypersec}

Recall that $H\sub A_{\th}$ is the set of elements $a$ such that $\tr(a)=0$.
For every $a\in A_{\th}$ consider the operator
$$d_a:=\de_{\tau}+\ad(a):A_{\th}\to A_{\th}.$$
We are interested in the following subset of $H$:
$$\Th=\Th_{\tau}=\{a\in H:\ \dim\ker(d_a)>1\}.$$
Note that we always have $\C\sub\ker(d_a)$, so $\Th$ coincides with
the set of $a\in H$ such that $\ker(d_a|_H)\neq 0$.
We are going to show that $\Th_{\tau}$ is a complex hypersurface in $H$
in the sense that it can be given locally as the zero set of a holomorphic
function on $H$. 

\begin{lem}\label{endlem} 
For every $a\in A_{\th}$ one has
$\ker(d_a)\simeq \Hom_{\CC}(E,E)$, $\coker(d_a)\simeq\Ext^1_{\CC}(E,E)$,
where $E$ is the holomorphic bundle $(A_{\th},\ov{\nabla}_a:=\de_{\tau}+a)$.
In particular, the operator $d_a$ is Fredholm of index zero.
The same assertions are true if we replace $d_a$ by its extension to
Sobolev spaces $d_a:W_{s+1}\to W_s$.
\end{lem}  

\Pf . The identification of $\ker(d_0)$ with endomorphisms of $E$
in the category $\CC$ follows from the definitions.
The identification of $\coker(d_0)$ with $\Ext^1_{\CC}(E,E)$ is 
constructed similarly to Proposition 2.4 of \cite{PS}). Note that
by Serre duality $\Ext^1_{\CC}(E,E)\simeq\Hom_{\CC}(E,E)^*$
(see \cite{P}).
The remaining assertions are proved using the same techniques
as in Theorem 2.8 of \cite{P}.
\ed

Let $\FF_0(H)$ be the space of Fredholm operators of index zero on 
$W_0(H)$, where $W_0(H)=\{a\in W_0(A_{\th})\ :\ \tr(a)=0\}$. 
Recall (see \cite{Qu},\cite{Fu})
that there is a natural holomorphic line bundle $\LL$ over 
$\FF_0(H)$ (called the {\it determinant line bundle}). By definition,
$\LL$ is trivial over each open set $\UU_A\sub \FF_0(H)$, where
$A:W_0(H)\to W_0(H)$ is a fixed trace class operator, 
$\UU_A$ consists of operators
$T$ such that $T+A$ is invertible. The transition functions of $\LL$
with respect to this open covering are given by the functions
$$g_{A,B}(T)=\operatorname{det}_F(1+(A-B)(T+B)^{-1})=
\operatorname{det}_F((T+A)(T+B)^{-1}),$$
where $A,B$ are trace class operators, and $\det_F$ denotes the Fredholm
determinant. 

Let $\ZZ\sub\FF_0(H)$ be the complement of $\UU_0$, i.e.,
the subset consisting of $T$ such that $\ker(T)\neq 0$.
Then $\ZZ$ is a complex hypersurface in $\FF_0(H)$.
More precisely, there exists a global holomorphic section $s$ of the
determinant line bundle $\LL$ such that $\ZZ$ coincides with the
zero locus of $s$. Under the standard trivialization of $\LL$ over
$\UU_A$ the section $s$ corresponds to the holomorphic function
$$s_A(T)=\operatorname{det}_F(1-(T+A)^{-1}A).$$
One can easily check that $s_B=g_{A,B}s_A$, so these functions glue
into a global section of $\LL$. One has $s_A(T)=0$ iff $1-(T+A)^{-1}A$
has a nonzero kernel. But $1-(T+A)^{-1}A=(T+A)^{-1}T$, so 
$s_A(T)=0$ iff $\ker(T)\neq 0$. This proves that the zero locus of $s$
coincides with $\ZZ$.

We have a map 
$$H\to\FF_0(H):a\mapsto d_a|_{W_0(H)}.$$
Note that the restriction of $d_a$ to $H$ (resp., $W_0(H)$) 
still has index zero, since 
$d_a(A_{\th})\sub H$ and $\ker(d_a)=\C\oplus\ker(d_a|_H)$.
Abusing the notation we denote by $\LL$ also 
the pull-back of the determinant line bundle to $H$ under the above map.

\begin{prop}\label{hyperprop} 
(i) The subset $\Th\sub H$ coincides with the zero locus
of the global holomorphic section $s(d_a|_{W_0(H)})$ of $\LL$ over $H$.

\noindent
(ii) For every one-dimensional subspace $L\sub H$ the intersection 
$\Th\cap L$ is a discrete subset of $L\setminus\{0\}$.
\end{prop}

\Pf .
(i) This is clear since $\Th$ is the preimage of $\ZZ\sub\FF_0(H)$.

\noindent
(ii)
We observe that the operator $\de_{\tau}:W_0(H)\to W_{-1}(H)$
is an isomorphism and the operator $\ad(a):W_0(H)\to W_{-1}(H)$
is compact for every $a\in A_{\th}$ (see the proof of Theorem 2.8 in \cite{P}).
Hence, for every $a\in H$ and $t\in\C$ we have $ta\in\Th$ iff
$-t^{-1}$ belongs to the spectrum of the compact operator
$\de_{\tau}^{-1}\ad(a):W_0(H)\to W_0(H)$.  
\ed

\begin{rem} Part (ii) of the above proposition implies that
$\Th$ contains no linear subspaces of $H$. On the other hand,
we will see later that $\Th$ is swept by infinite-dimensional
{\it affine} subspaces (see Theorem \ref{projthm}).
\end{rem} 

We can also consider more general loci
$$\Th^{(n)}_{\tau}=
\{a\in H:\ \dim\ker(d_a)>n+1\}=
\{a\in H:\ \dim\ker(d_a|_H)>n\}$$
for $n\ge 0$. For $n=0$ we get $\Th^{(0)}_{\tau}=\Th_{\tau}$.
Since the kernel of $d_a$ coincides with the kernel of its extension
to $W_0$, it follows that the loci $\Th^{(n)}_{\tau}$ are closed.

\subsection{Logarithmic  derivative and connected components of groups
of invertible elements}

In this section we will study the relation between holomorphic structures and 
the 
group of invertible elements in $A_{\th}$.

\begin{thm}\label{logthm} Consider the map
$$L_{\tau}: A_{\th}^*\to A_{\th}: x\mapsto x^{-1}\de_{\tau}(x).$$

\noindent
(i) The composition
$$\chi_{\tau}:=\tr\circ L_{\tau}:A_{\th}^*\to\C$$ 
is a locally constant homomorphism with image $(2\pi i)(\Z+\Z\tau)\sub\C$
and kernel
$(A_{\th}^*)_0$, the connected component of $1$ in $A_{\th}^*$.

\noindent
(ii) Let $H\sub A_{\th}$ be the set of elements $a$ such that $\tr(a)=0$.
Let 
$$\Om_{\tau}=\{ a\in H:\ \dim\ker(\de_{\tau}+\ad(a))=1\}$$
be the complement to the hypersurface $\Th_{\tau}\sub H$. 
Then $\Om_{\tau}$ is a dense open subset of $H$ and
$L_{\tau}$ induces homeomorphisms
$$A_{\th}^*/\C^*\wt{\to}\Om_{\tau}+2\pi i (\Z+\Z\tau),$$
$$(A_{\th}^*)_0/\C^*\wt{\to}\Om_{\tau}.$$  
\end{thm}

\Pf . From the Leibnitz rule we immediately get
\begin{equation}\label{Leq}
L_{\tau}(xy)=y^{-1}L_{\tau}(x)y+L_{\tau}(y).
\end{equation}
Taking traces we see that $\chi_{\tau}=\tr\circ L_{\tau}$ is a homomorphism.
It is easy to see that the derivative of $L_{\tau}$ at the point
$x\in A_{\th}^*$ is $a\mapsto x^{-1}\de_{\tau}(ax^{-1})x$.
Hence, the derivative of $\chi_{\tau}$ is
$a\mapsto \tr(x^{-1}\de_{\tau}(ax^{-1})x)=0$.
This implies that $\chi_{\tau}$ is locally constant. In particular,
$\chi_{\tau}((A_{\th}^*)_0)=0$, i.e., $L_{\tau}((A_{\th})_0)\in H$.

Next we claim that the map $L_{\tau}:(A_{\th}^*)_0/\C^*\to H$ is a local
homeomorphism. Indeed, (\ref{Leq}) shows that it is enough to
check that $L_{\tau}$ is a local homeomorphism in a neighborhood of
$1\in (A_{\th}^*)_0/\C^*$. Then we can use the map $\bfe_l$ constructed in
Theorem \ref{convthm} to get a local inverse to $L_{\tau}$.

We know that the space of holomorphic
endomorphisms of $(A_{\th},\de_{\tau}+a)$
can be identified
with $\ker(\de_{\tau}+\ad(a))$ (see Lemma \ref{endlem}). 
Using the classification of holomorphic
bundles on noncommutative tori given in \cite{P} we derive that
the dimension of this space is $1$ iff there exists a holomorphic
isomorphism of $(A_{\th},\ov{\nabla})$ with a standard holomorphic
bundle $(A_{\th},\de_{\tau}+z)$ for some $z\in\C$. In other words,
$\ker(\de_{\tau}+\ad(a))=\C$ iff there exists $x\in A_{\th}^*$ such
that $a\equiv L_{\tau}(x)\mod\C$.
Thus, the map 
\begin{equation}\label{mapeq}
A_{\th}^*/\C^*\to H:x\mapsto L_{\tau}(x)-\tr(L_{\tau}(x))
\end{equation}
has $\Om_{\tau}$ as an image. As we have seen above this map is
a local homeomorphism, so $\Om_{\tau}$ is open in $H$.
Now we claim that nonempty fibers of (\ref{mapeq})
are exactly orbits of the action of the central subgroup
$\Z^2\sub A_{\th}^*/\C^*$. Indeed,
assume that $L_{\tau}(x)=L_{\tau}(y)\mod\C$ for some $x,y\in A_{\th}^*$.
Set $z=xy^{-1}$. Then
$$L_{\tau}(x)=L_{\tau}(zy)=y^{-1}L_{\tau}(z)y+L_{\tau}(y)$$
which implies that $y^{-1}L_{\tau}(z)y\in\C$, hence, $L_{\tau}(z)\in\C$.
But this is possible only if $z$ is proportional to $U_v$ for some
$v\in\Z^2$. This proves our claim about the fibers of the map (\ref{mapeq}).

Proposition \ref{hyperprop}(ii) immediately implies that 
$\Om_{\tau}$ is connected (and dense in $H$). Since
we have identified the quotient of $A_{\th}^*/\C^*$ by $\Z^2$ with 
$\Om_{\tau}$, it follows that the embedding $\Z^2\sub A_{\th}^*$
induces an isomorphism on connected components. Since $\chi_{\tau}$
is locally constant this implies that 
$$\chi_{\tau}(A_{\th}^*)=\chi_{\tau}(\Z^2)=(2\pi i)(\Z+\Z\tau).$$
We also see that $L_{\tau}:(A_{\th}^*)_0/\C^*\to\Om_{\tau}$ 
is a homeomorphism.
\ed

\begin{cor}\label{maineqcor} 
For every $a\in H$ there exists a nonzero $x\in A_{\th}$ such that
$\de_{\tau}(x)=ax$ (resp., $\de_{\tau}(x)=xa$).
\end{cor}

\Pf . If $-a\in\Om_{\tau}$ then by the above theorem there exists $x\in 
A_{\th}^*$
such that $\de_{\tau}(x)=-xa$. Therefore, $\de_{\tau}(x^{-1})=ax^{-1}$.
If $-a\not\in\Om_{\tau}$ 
then the holomorphic bundle $(A_{\th},\ov{\nabla}_{-a}=\de_{\tau}-a)$
is decomposable. Since one of the indecomposable factors should be of positive 
degree 
this implies that $H^0(A_{\th},\ov{\nabla}_{-a})\neq 0$.
The case of the equation $\de_{\tau}(x)=xa$ reduces to the previous
case by using the identification $A_{\th}^{opp}=A_{-\th}$.
\ed

As another corollary we get a new proof of the following well known fact.

\begin{cor} The embedding $\Z^2\sub A_{\th}^*$ induces an isomorphism
on connected components. 
\end{cor}

Using the results of \cite{P} we can give one more characterization of 
the open subset $\Om_{\tau}\sub H$. As in section \ref{hypersec} this
leads to a set of local equations for its complement $\Th_{\tau}$.
For every $a\in H$ let us set
$\ov{\nabla}_a=\de_{\tau}+a$.

\begin{prop} For $a\in H$ one has $a\in\Om_{\tau}$ iff there exists
$z\in\C$ such that 
$$\ker(\ov{\nabla}_{a+z})=\{x:\ \de_{\tau}(x)+(a+z)x=0\}$$
is zero. Moreover, for $a\in\Om_{\tau}$ this holds for
all $z\not\in 2\pi i(\Z+\Z\tau)$.
Hence, if we choose any $z_0\not\in 2\pi i(\Z+\Z\tau)$ then
$a\in\Om_{\tau}$ iff the operator $\ov{\nabla}_{a+z_0}$ is invertible.
\end{prop}

\Pf . If $a\in \Om_{\tau}$ then this follows from the fact that cohomology of
the standard holomorphic bundle $(A_{\th},\de_{\tau}+z)$ vanishes for 
$z\not\in 2\pi i(\Z+\Z\tau)$.
If $a\not\in\Om_{\tau}$ then it suffices to observe that the holomorphic bundle
$(A_{\th},\ov{\nabla}_{a+z})$ has a direct factor of positive degree.
\ed

One can generalize some of the assertions of Theorem \ref{logthm} to
a slightly more general context.

\begin{prop}\label{holbunprop} 
Let $(E,\ov{\nabla})$ be a basic right module over $A_{\th}$ equipped
with a holomorphic structure. Then the map
$$\chi_E:\Aut_{A_{\th}}(E)\to\C:x\mapsto \tr(x^{-1}[\ov{\nabla},x])$$ 
is a locally constant homomorphism that does not depend on $\ov{\nabla}$.
Its kernel coincides with the connected component of $1$ in $\Aut_{A_{\th}}(E)$
and its image is the lattice $\frac{2\pi i}{\rk(E)}(\Z+\tau\Z)\sub\C$, where
$\rk(E)\in(\Z+\Z\th)\cap \R_{>0}$ is the rank of $E$.
If $E=E'\oplus E''$ is a decomposition into the direct sum of $A_{\th}$-modules, 
where $E'$ is also basic, then
one has the commutative diagram
\begin{equation}\label{comdiag}
\begin{diagram}
\Aut_{A_{\th}}(E') &\rTo{\chi_{E'}} &\C \\
\dTo &&\dTo_{\rk(E')/\rk(E)}\\
\Aut_{A_{\th}}(E) &\rTo{\chi_E} &\C
\end{diagram}
\end{equation}
\end{prop}

\Pf . Any other holomorphic structure on $E$ has form $\ov{\nabla}+\phi$
for some $\phi\in\End_{A_{\th}}(E)$. But
$$\tr(x^{-1}[\ov{\nabla}+\phi,x])=\tr(x^{-1}[\ov{\nabla},x])+
\tr(x^{-1}\phi x-\phi)=\tr(x^{-1}[\ov{\nabla},x]).$$
This shows that $\chi_E$ is independent of $\ov{\nabla}$. If we choose
$\ov{\nabla}$ to be standard (see \cite{P}) then the pair
$(\End_{A_{\th}},\ad\ov{\nabla})$ can be identified with
a pair of the form $(A_{\th'},\de_{\tau}/r)$ for some $\th'$, where
$r=\rk(E)$ (see \cite{PS}, Prop.~2.1). It remains to apply Theorem \ref{logthm}.

Since the homomorphism $\chi_E$ does not depend on a choice of a holomorphic 
structure
we can choose a holomorphic structure on $E$ compatible with the decomposition
$E=E'\oplus E''$. Then commutativity of (\ref{comdiag}) follows from the 
compatibility
of the embedding $i:\End_{A_{\th}}(E')\to\End_{A_{\th}}(E)$ with normalized 
traces:
$$\tr(i(x))=\frac{\rk(E')}{\rk(E)}\tr(x).$$
\ed

As a corollary we get a new proof of the following well known fact.

\begin{cor}\label{isomcor} 
Let $E$ be a basic right module over $A_{\th}$, and let
$E=E'\oplus E''$ be a decomposition into the direct sum of $A_{\th}$-modules, 
where $E'$ is also a basic module.
Then the natural homomorphism $\Aut_{A_{\th}}(E')\to\Aut_{A_{\th}}(E)$
induces an isomorphism on connected components.
\end{cor}

\Pf . Apply commutative diagram \eqref{comdiag} together with the fact that
the image of $\chi_E$ (resp., of $\chi_{E'}$) is $2\pi i\rk(E)(\Z+\Z\tau)$
(resp., $2\pi i\rk(E')(\Z+\Z\tau)$).
\ed

\begin{cor}\label{surcor} 
Let $E$ be a basic right module over $A_{\th}$, and let $E=E'\oplus E''$
be a decomposition into a direct sum of $A_{\th}$-modules. Then the natural
homomorphism $\Aut_{A_{\th}}(E')\to\Aut_{A_{\th}}(E)$ induces a surjection on
connected components. 
\end{cor}

\Pf . It suffices to choose any decomposition $E'=E'_1\oplus E'_2$ with $E'_1$ 
basic
and apply Corollary \ref{isomcor}.
\ed

\begin{rem} The well known fact that the embedding of $A_{\th}$ into 
$\Mat_n(A_{\th})$
induces an isomorphism on connected components of 
groups of invertible elements (for $\ov{A}_{\th}$ this statement is Theorem
8.3 of \cite{R-canc}, the case of $A_{\th}$ follows 
using its invariance under the holomorphic functional calculus)
allows to deduce from Corollary \ref{isomcor} 
the same statement for $E'$ not necessarily basic.
\end{rem}

\subsection{More on discriminant hypersurface}

It would be interesting to study intersections of $\Th_{\tau}$ with
finite-dimensional subspaces of $H$. In the following proposition
we consider simplest examples.

\begin{prop}\label{linesprop}
(i) Let $h:\Z^2\to\R$ be a homomorphism. Assume that
$a\in H$ is an element such that $h>\eps$ on $\supp(a)$
for some $\eps>0$. Then $\C a\in\Om_{\tau}$. The same conclusion
holds if we assume that $h$ takes values in $\Q$ and that
$h\ge 0$ on $\supp(a)$.

\noindent
(ii)
Let $e\in A_{\th}$ be a nontrivial idempotent such that $\de_{\tau}$ preserves
$eA_{\th}\sub A_{\th}$, i.e., $e\de_{\tau}(e)=\de_{\tau}(e)$.
Then $\Th_{\tau}\cap \C\de_{\tau}(e)=\{\de_{\tau}(e)\}$.
\end{prop}

We need the following simple lemma for the proof.

\begin{lem}\label{uniquelem} 
Let $a=x_0^{-1}\de_{\tau}(x_0)$ for some $x_0\in A_{\th}^*$.
Then every $x\in A_{\th}$ such that $\de_{\tau}(x)=xa$
is proportional to $x_0$.
\end{lem}

\Pf . Indeed, we have
$$\de(xx_0^{-1})=\de(x)x_0^{-1}-xx_0^{-1}\de(x_0)^{-1}x_0^{-1}=0,$$
hence $xx_0^{-1}$ is a constant.
\ed 

\noindent
{\it Proof of Proposition \ref{linesprop}}. 
(i) It suffices to show that the equation $L_{\tau}(x)=a$ has a solution
with $x\in A_{\th}^*$.
If $h>\eps$ on $\supp(a)$ then this follows
from Proposition \ref{suppconvprop}.
If $h$ takes values in $\Q$ and $h\ge 0$ on $\supp(a)$ then we can write
$a=b+c$, where $h=0$ on $\supp(b)$ and $h>\eps>0$ on $\supp(c)$. Since 
the subalgebra of elements supported on $h=0$ is commutative and stable
under $\de_{\tau}$, we have $b=L_{\tau}(x)$ for some $x\in A_{\th}^*$ such that
$h=0$ on $\supp(x)$. Now let us consider the element $c'=xcx^{-1}$.
Note that $h>0$ on $\supp(c')$. Therefore, by Proposition \ref{suppconvprop}
there exists an element $y\in A_{\th}^*$ such that $c'=L_{\tau}(y)$.
Then
$$L_{\tau}(yx)=x^{-1}L_{\tau}(y)x+L_{\tau}(x)=x^{-1}c'x+b=c+b=a.$$

\noindent
(ii)
For $\la\in\C^*$ consider the element $x_{\la}=\la e+(1-e)\in A_{\th}$.
Then $x_{\la}\in (A_{\th}^*)_0$ and $x_{\la}^{-1}=\la^{-1}e+(1-e)$. Hence,
$$L_{\tau}(x_{\la})=
(\la^{-1}e+(1-e))(\la-1)\de_{\tau}(e)=(1-\la^{-1})\de_{\tau}(e)
$$
It follows that $z\de_{\tau}(e)\in\Om_{\tau}$ for all $z\in\C\setminus\{1\}$.
Now let us prove that $\de_{\tau}(e)\not\in\Om_{\tau}$.
Assume that there exists $x_0\in A_{\th}^*$ such that 
$$x_0^{-1}\de_{\tau}(x_0)=\de_{\tau}(e).$$
Since we also have $\de_{\tau}(e)=e\de_{\tau}(e)$, applying Lemma 
\ref{uniquelem} to $x=e$ we derive that $e$ should be proportional
to $x_0$. Since $e$ is not invertible we get a contradiction.
\ed

\begin{rem} There are plenty of idempotents $e\in A_{\th}$ such that
$eA_{\th}$ is preserved by $\de_{\tau}$. In fact,
for every number $0<m\th+n<1$ such that $m, n\in\Z$ and $m<0$ there exists
an idempotent $e$ as above with $\tr(e)=m\th+n$ (see \cite{P-holomid}).
\end{rem}

Next, we are going to describe the decomposition discriminant locus
$\Th_{\tau}\sub H$ corresponding to types of idempotents in $A_{\th}$.

\begin{lem} Assume that $\th$ is irrational. Then two idempotents
$e$ and $e'$ in $A_{\th}$ are conjugate by an element in $(A_{\th}^*)_0$ iff 
$\tr(e)=\tr(e')$.
\end{lem}

\Pf . The ``only if'' part is trivial. Assume that $\tr(e)=\tr(e')$.
Then by Rieffel's classification of projective $A_{\th}$-modules
(see \cite{R-canc}) there exist isomorphisms of right $A$-modules
$eA\simeq e'A$ and $(1-e)A\simeq (1-e')A$. Therefore, there exists
an element $x\in A_{\th}^*$ such that $xeA=e'A$ and $x(1-e)A=(1-e')A$.
This immediately implies the equality of the idempotents 
$xex^{-1}=e'$. It remains to show that $x$ can be chosen in $(A_{\th}^*)_0$.
But this follows easily from Corollary \ref{surcor}.
\ed

For every idempotent $e$ 
let us denote by $A_{e,\tau}$ the set of $a\in A_{\th}$
such that the operator $x\mapsto \de(x)+ax$ preserves the decomposition
$A_{\th}=eA_{\th}\oplus (1-e)A_{\th}$, i.e.,
$$A_{e,\tau}=\{a\in A_{\th}|\ \de(e)+ae\in eA_{\th}, 
\de(1-e)+a(1-e)\in (1-e)A_{\th}\}.$$
Equivalently, $a\in A_{e,\tau}$ iff $a$ satisfies the equations
$$e\de(e)=ea(1-e), \ (1-e)ae=(1-e)\de(1-e).$$
This shows that $A_{e,\tau}$ is an affine subspace in $A_{\th}$ with
the associated linear subspace $eAe+(1-e)A(1-e)$.

\begin{lem} For $x\in A_{\th}^*$ and for an idempotent $e\in A_{\th}$
the map $a\mapsto x^{-1}ax+L_{\tau}(x)$ induces a bijection from
$A_{e,\tau}$ to $A_{x^{-1}ex,\tau}$.
\end{lem}

\Pf .
In fact, the natural action of $A_{\th}^*$ on operators of the
form $\de_{\tau}+a$ is equivalent to the following action of
$A_{\th}^*$ on $A_{\th}$: 
$$x*a=xax^{-1}+L_{\tau}(x^{-1}),$$
where $x\in A_{\th}^*$, $a\in A_{\th}$. 
It is clear from the definition that the action of $x^{-1}$ sends
$A_{e,\tau}$ to $A_{x^{-1}ex,\tau}$.
\ed

We refer to the action of $A_{\th}^*$ on $A_{\th}$ introduced above
as {\it twisted action}. Note that it preserves $H$ and $\Th_{\tau}$.

Let us set $H_{e,\tau}=A_{e,\tau}\cap H$.
For every number $r\in(0,1)\cap(\Z+\Z\th)$ let us set
$$\HH_{r,\tau}=\cup_{e:\tr(e)=r} H_{e,\tau},$$
where the union is taken over the set of all idempotents $e$ with $\tr(e)=r$.
Since $A_{1-e,\tau}=A_{e,\tau}$ we have $\HH_{1-r,\tau}=\HH_{r,\tau}$.

\begin{thm}\label{projthm}
Assume that $\th$ is irrational. Then
\begin{equation}\label{uneq}
\Th_{\tau}=\cup_{r\in (0,1/2)\cap(\Z+\Z\th)}\HH_{r,\tau}.
\end{equation}
This decomposition is irreducible, i.e., none of the subsets
is contained in the union of the rest.
Also if $e\in A_{\th}$ is any idempotent with $\tr(e)=r$ then
we have a map
\begin{equation}\label{Admap}
\psi_e:(A_{\th}^*)_0\times H_{e,\tau}\to :(x,a)\mapsto x^{-1}ax+L(x)
\end{equation}
such that $\HH_{r,\tau}$ is the image of $\psi_e$.
The codimension of the differential of $\psi_e$ at a point
$(x,a)$ is equal to $\dim\Hom_{\CC}(E_1,E_2)+\dim\Hom_{\CC}(E_2,E_1)$,
where $\CC$ is the category of holomorphic bundles on $T_{\th,\tau}$,
$E_1$ and $E_2$ are the summands of the holomorphic decomposition
$(A_{\th},\de_{\tau}+\psi(x,a))=E_1\oplus E_2$ corresponding to the
idempotent $x^{-1}ex$. 
\end{thm}

\Pf . From the main theorem in \cite{P} we know that
a holomorphic bundle $E$ on a noncommutative torus has non-scalar
endomorphisms then it is either decomposable or its rank is a non-primitive
element of $\Z+\Z\th$. Since for a holomorphic bundle 
$E_a=(A_{\th},\de_{\tau}+a)$
the second alternative is impossible, we conclude that $a\in\Th_{\tau}$
iff $E_a$ is decomposable. This implies \eqref{uneq}.
Since $\psi_e$ is $(A_{\th}^*)_0$-equivariant with respect to the
twisted action on $H$ and the action by left multiplication on
$(A_{\th}^*)_0$, it suffices to consider the differential
of $\psi_e$ at a point $(1,a_0)$, where $a_0\in H_e$.
This differential maps $(x,a)\in A_{\th}\oplus H\cap(eAe+(1-e)A(1-e))$ to
$\de_{\tau}(x)+[a_0,x]+a=d_{a_0}(x)+a$.
But by Lemma \ref{endlem}
$$\coker(d_{a_0}:A_{\th}\to A_{\th})\simeq\Ext^1_{\CC}(E,E),$$
where $E=(A_{\th},\de_{\tau}+a_0)$. 
Furthermore, it is easy to see that the map $d_{a_0}$ is 
compatible with the decomposition 
$$A_{\th}=eA_{\th}e+(1-e)A_{\th}(1-e)\oplus eA_{\th}(1-e)\oplus (1-e)A_{\th}e$$
and that the cokernel of $d_{a_0}$ acting 
on each of these four pieces is identified with $\Ext^1(E_i,E_j)$
where $E=E_1\oplus E_2$ is the holomorphic decomposition induced by $e$.
The assertion immediately follows from this.
\ed

Let $r_0$ be the unique number of the form $\pm\th+n$ with $n\in\Z$ such
that $0<r_0<1/2$. 

\begin{prop} One has 

\noindent
(i) $\Th_{\tau}^{(1)}=\Th_{\tau}$;

\noindent
(ii) $\Th_{\tau}^{(2)}=
\cup_{r\in (0,1/2)\cap(\Z+\Z\th), r\neq r_0}\HH_{r,\tau}$;

\noindent
(iii) the natural map 
$$\sqcup_{e:\tr(e)=r_0}
H_{e,\tau}\setminus\Th_{\tau}^{(2)}\to \Th_{\tau}\setminus\Th_{\tau}^{(2)}$$
is an affine line bundle.
\end{prop}

\Pf . (i) follows from the fact that if $E_a$ is decomposable
then $\dim\End_{\CC}(E_a)>2$. Moreover, the only case when
$\dim\End_{\CC}(E_a)=3$ is the following: $E_a\simeq E_1\oplus E_2$
in the holomorphic category,
where $E_1$ and $E_2$ are simple, $\Hom_{\CC}(E_2,E_1)=0$ and 
$\Hom_{\CC}(E_1,E_2)$ is one-dimensional. Let $\rk E_1=m\th+n$.
Then $\rk E_2=-m\th+(1-n)$ and hence
$$\dim\Hom_{\CC}(E_1,E_2)-\dim\Hom_{\CC}(E_2,E_1)=-mn-m(1-n)=-m$$
(see Corollary 2.9 of \cite{P}). 
This implies that $m=-1$, so either $\rk E_1=r_0$ or 
$\rk E_2=r_0$ which proves (ii). It remains to note that the set
of all possible holomorphic
decomposition of $E_a$ into $E_1$ and $E_2$ is
a principal homogeneous space for $\Hom_{\CC}(E_1,E_2)$ 
which leads to (iii).
\ed

The part of Theorem \ref{projthm} concerning tangent maps
should imply that $\HH_{r,\tau}$ has (in an appropriate
sense) codimension $\ge 2$ for
$r\neq r_0$ (where $0<r<1/2$)
while $\HH_{r_0,\tau}$ is irreducible of codimension $1$. If we were
in the finite-dimensional situation we would immediately deduce from this
that the hypersurface $\Th_{\tau}$ is irreducible.
Our lack of knowledge does not allow us to state this precisely. 
However, we conjecture that the above statements hold
when we intersect these loci with generic finite-dimensional subspaces in
$A_{\th}$.

\end{document}